\documentstyle[leqno]{article}
\makeindex
\title{Complex almost contact metric structures on complex hypersurfaces in hyperk\"{a}hler manifolds}
\author{MITSUHIRO IMADA 
\\ National Institute of Technology, Ibaraki Collage}
\begin{document}
\renewcommand{\thefootnote}{\fnsymbol{footnote}}
\maketitle
Abstract: In this paper, we construct complex almost contact metric structures on complex hypersurfaces in hyperk\"{a}hler manifolds. This construction is analogous to that in contact geometry. \\ \\

\section{Introduction}
\ \ The theory of complex contact geometry started with the papers of Kobayashi $\cite{Kob}$, as a variant of real contact geometry.  More recent examples, including complex projective space and the complex Heisenberg group, are given in $\cite{B}$. Ishihara and Konishi $\cite{IK2}$ defined the so-called I-K normality of complex contact manifolds as for Sasakian manifolds in real contact geometry. In this paper, we construct complex almost contact manifolds from hyperk\"{a}hler manifolds. Leaving the detailed notion of hyperk\"{a}hler manifolds to Definition 3.1, we state the first main result as follows: \\ \\ 
$\bf{Theorem \ A} \ (Theorem \ 3.3.)$ \  Let $(\widetilde{M}, J_1, J_2, J_3, \tilde{g})$ be a hyperk\"{a}hler manifold and $M$ be a complex hypersurface of $\widetilde{M}$. The inclusion $\iota : M \longrightarrow \widetilde{M}$ canonically induces a complex almost contact metric structure on $M$. \\ \\
 \ \ This main result is an analogous to Morimoto $\cite{Mo}$. He shows that real hypersurfaces in K\"{a}hler manifolds equip an almost contact metric structure induced from the ambient K\"{a}hler structure. \\ 
 \ \ In addition, we show that covariant derivatives of tensors belonging to complex almost contact metric structures have the following forms.  \\ \\
$\bf{Theorem \ B} \ (Theorem \ 4.1.)$ \  Let $(G, H, J, u, v, U, V, g)$ be a complex almost contact metric structure on a complex hypersurface $M$. Then the derivatives of $G$ and $H$ have the following forms: 
\begin{eqnarray*}
&{ }& (\nabla_X G) Y = - u(Y)AX + v(Y)JAX + g(AX, Y)U - g(JAX, Y)V, \\
&{ }& (\nabla_X H) Y = - u(Y)JAX - v(Y)AX + g(AX, Y)V + g(JAX, Y)U. 
\end{eqnarray*}
 \ \ This paper is organized as follows. In section 2, we recall definitions of complex contact manifolds and hyperk\"{a}hler manifolds. In section 3, we prove the main theorem which constructs complex almost contact metric structures on complex hypersurfaces in hyperk\"a{h}ler manifolds. In section 4, we give some results of tensor calculations of these complex almost contact metric structures. 
\section{Definitions}
\ \ \ \ \ We first recall the notion of complex contact metric manifolds $\cite{B}$. \\
$\bf{Definition \ 2.1. \ }$ Let $M$ be a complex manifold with dim$_{\mathbf{C}} M =2n+1$ and $J$ the complex structure on $M$.  $M$ is called a complex contact manifold if there exists an open covering $\mathcal{U}=\{ \mathcal{O}_{\lambda} \}$ of $M$ such that: \\ \\
\ \ \ \ \ \ \ 1) On each $\mathcal{O}_{\lambda}$ there is a holomorphic 1-form $\omega_{\lambda}$ with $\omega_{\lambda} \wedge (d\omega_{\lambda})^n \neq 0 \ $ everywhere; \\
\ \ \ \ \ \ \ 2) If $\mathcal{O}_{\lambda} \cap \mathcal{O}_{\mu} \neq \phi $, there is a nonvanishing holomorphic function $h_{\lambda \mu }$ on $\mathcal{O}_{\lambda} \cap \mathcal{O}_{\mu}$ such that
\begin{eqnarray}
 \omega_{\lambda} = h_{\lambda \mu }\omega_{\mu} \ \ {\rm in} \ \mathcal{O}_{\lambda} \cap \mathcal{O}_{\mu}.
\end{eqnarray} 

For each $\mathcal{O}_{\lambda}$, we define a distribution $\mathcal{H}_{\lambda}$ = $\{ X \in T\mathcal{O}_{\lambda}$  $|  \ \omega_{\lambda}(X)=0 \}$. Note that the $h_{\lambda \mu}$ are nonvanishing, and $\mathcal{H}_{\lambda} = \mathcal{H}_{\mu}$ on $\mathcal{O}_{\lambda} \cap \mathcal{O}_{\mu}$. Thus $\mathcal{H} = \cup \mathcal{H}_{\lambda}$ is a holomorphic, nonintegrable subbundle on $M$, called the horizontal subbundle. \\ \\ 
$\bf{Definition \ 2.2. } $ Let $M$ be a complex manifold with dim$_{\mathbf{C}}=2n+1$ and $J$ a complex structure. Let $g$ be a Hermitian metric. $M$ is called a complex almost contact metric manifold if there exists an open covering $\mathcal{U}=\{ \mathcal{O}_{\lambda} \}$ of $M$ such that:  \\ \\
1)  On each $\mathcal{O}_{\lambda}$ there are 1-forms $u_{\lambda}$ and $v_{\lambda}= - u_{\lambda}J$, (1,1) tensors $G_{\lambda}$ and $H_{\lambda}= JG_{\lambda}$, unit vector fields $U_{\lambda}$ and $V_{\lambda}= J U_{\lambda}$ such that 
\begin{eqnarray}
&{ }& G_{\lambda} J_{\lambda}=-J_{\lambda} G_{\lambda}, \ \ \ H_{\lambda}^2=G_{\lambda}^2= -id + u_{\lambda} \otimes U_{\lambda} + v_{\lambda} \otimes V_{\lambda}, \nonumber \\
&{ }& g(G_{\lambda}X, Y) = -g(X, G_{\lambda} Y), \ \ \ g(U_{\lambda}, X)=u_{\lambda}(X),  \\
&{ }& G_{\lambda} U_{\lambda}=0, \ \ \ u_{\lambda} (U_{\lambda})=1; \nonumber
\end{eqnarray}
2) If $\mathcal{O}_{\lambda} \cap \mathcal{O}_{\mu} \neq \phi $, there are functions $a,b$ on $\mathcal{O}_{\lambda} \cap \mathcal{O}_{\mu}$ such that

\begin{eqnarray}
&{ }& u_{\mu}=au_{\lambda} - b v_{\lambda}, \ \ \ v_{\mu}=bu_{\lambda} + av_{\lambda}, \nonumber  \\
&{ }& G_{\mu}=aG_{\lambda} - bH_{\lambda}, \ \ \ H_{\mu}=bG_{\lambda} + aH_{\lambda}, \\
&{ }& a^2+b^2=1.  \nonumber 
\end{eqnarray} 
$\bf{Definition \ 2.3.}$ Let $(M, \{ \omega _{\lambda}\})$ be a complex contact manifold with complex contact structure $J$ and Hermitian metric $g$. We call $(M, J, G, u, U, g)$ a complex contact metric manifold if there exists an open covering $\mathcal{U}=\{ \mathcal{O}_{\lambda} \}$ of $M$ such that (here and below $G=G_{\lambda},$ etc) :  \\ \\
\ \ \ 1) On each $\mathcal{O}_{\lambda}$ there is a local (1,1) tensor $G_{\lambda}$ such that $(u_{\lambda}, v_{\lambda}, U_{\lambda}, V_{\lambda}, G_{\lambda}, H_{\lambda}=G_{\lambda}J, g)$ is an almost contact metric structure on $M$; \\ \\
\ \ \ 2) $g(X, G_{\lambda}Y) = du_{\lambda}(X, Y) + (\sigma_{\lambda} \wedge v_{\lambda}) (X,Y)$ and $g(X, H_{\lambda}Y) = dv_{\lambda}(X, Y) - (\sigma_{\lambda} \wedge u_{\lambda}) (X,Y),$ where $\sigma_{\lambda} (X) = g(\nabla _{X}U_{\lambda}, V_{\lambda})$ with $\nabla$ the Levi-Civita connection with respect to $g$. \\ \\
$\bf{Remark \ 2.4.}$ Foreman $\cite{Fo}$ showed the existence of complex contact metric structures on complex contact manifolds. \\ \\ 
$\bf{Remark \ 2.5.}$ We can locally choose orthonormal vectors $X_{1}, \cdots , X_{n}$ in $\mathcal{H}$ such that $ \{ X_{i}, JX_{i}, GX_{i}, HX_{i}, U, V \ | \ 1 \leq i \leq n \} $ is an orthonormal basis of the tangent spaces of $U_{\alpha}$. \\ \\
We recall the definition of I-K normality introduced by Ishihara and Konishi $\cite{IK1}$ for (almost) complex contact metric structures. We set the two tensor fields $S$ and $T$ by,
\begin{eqnarray}
& S(X,Y)  &  = \ [G,G](X,Y) + 2g(X,GY)U - 2g(X,HY)V  \\
& { } & \ \ \ \ + 2v(Y)HX - 2v(X)HY + \sigma (GY)HX  \nonumber \\
& { } & \ \ \ \ - \sigma (GX)HY + \sigma (X)GHY - \sigma (Y)GHX, \nonumber \\
& T(X,Y)  &  = \ [H,H](X,Y) - 2g(X,GY)U +2g(X,HY)V  \\
& { } &  \ \ \ \ + 2u(Y)GX - 2u(X)GY + \sigma (HX)GY \nonumber \\
& { } &  \ \ \ \ - \sigma (HY)GX + \sigma (X)GHY - \sigma (Y)GHX. \nonumber
 \end{eqnarray} 
$\bf{Definition \ 2.6.}\ $ A complex contact manifold $M$ is I-K normal if the tensors $S$ and $T$ both  vanish. \\ \\
$\bf{Remark \ 2.7.} \ $ I-K normality implies that the underlying Hermitian manifold $(M, J, g)$ is a K\"{a}hler manifold (cf. $\cite{IK2}$). \\ \\

\section{Constructions}

\ \ In this chapter, we construct a complex almost contact metric structure on complex hypersurfaces in hyperk\"{a}hler manifolds. At first, we recall the definition of hyperk\"{a}hler manifolds. \\ \\
$\bf{Definition \ 3.1.}$ \  $(M, J_1, J_2, J_3, g)$ is a hyperk\"{a}hler manifold if $J_1, J_2, J_3$ are complex structures on a complex manifold $M$ satisfying 
\begin{eqnarray*}
J^2_1 = J^2_2 = J^2_3 =  J_1J_2J_3 = -id,
\end{eqnarray*}
$\bf{Definition \ 3.2.}\ $ A complex contact manifold $M$ is I-K normal if the tensors $S$ and $T$ both  vanish. \\ \\
Let $(\widetilde{M}, J_1, J_2, J_3, \tilde{g})$ be a hyperk\"{a}hler manifold and $M$ be a complex hypersurface in $\widetilde{M}$. For each $\widetilde{X} \in T\widetilde{M}$, we decompose $\widetilde{X}$ as follows:
\begin{eqnarray}
\widetilde{X} = X + \tilde{g} (\widetilde{X}, \xi) \xi + \tilde{g} (\widetilde{X}, J_1\xi) J_1\xi
\end{eqnarray}
where $X$ is the component of $\widetilde{X}$ tangent to $M$, and $\xi$ and $J_1\xi$ are normal to $M$. Applying $J_1$ to $(6)$, we have
\begin{eqnarray*}
J_1\widetilde{X} = J_1X + \tilde{g} (\widetilde{X}, \xi) J_1\xi - \tilde{g} (\widetilde{X}, J_1\xi) \xi.
\end{eqnarray*}
Applying $J_2$ to $(6)$, we have
\begin{eqnarray*}
J_2\widetilde{X} = J_2X + \tilde{g} (\widetilde{X}, \xi) J_2\xi - \tilde{g} (\widetilde{X}, J_1\xi) J_3\xi.
\end{eqnarray*}
Since $J_2\xi$ is not tangent to $TM$, we decompose $J_2X$ as follows:
\begin{eqnarray}
J_2 X = GX + \tilde{g} (J_2 X, \xi) \xi + \tilde{g} (J_2 X, J_1\xi) J_1\xi,
\end{eqnarray}
where $GX$ is the component of $J_2 X$ tangent to $TM$. Then
\begin{eqnarray}
J_2\widetilde{X} & = & GX + \tilde{g} (J_2 X, \xi) \xi + \tilde{g} (J_2 X, J_1\xi) J_1\xi + \tilde{g} (\widetilde{X}, \xi) J_2\xi - \tilde{g} (\widetilde{X}, J_1\xi) J_3\xi \nonumber \\
& = & GX + \tilde{g} (\widetilde{X}, \xi) J_2\xi - \tilde{g} (\widetilde{X}, J_1\xi) J_3\xi - \tilde{g} (X, J_2\xi) \xi + \tilde{g} (X, J_3\xi) J_1\xi \nonumber \\
& = & GX + \tilde{g} (\widetilde{X}, \xi) J_2\xi - \tilde{g} (\widetilde{X}, J_1\xi) J_3\xi - g(X, J_2\xi) \xi + g(X, J_3\xi) J_1\xi. 
\end{eqnarray}
Again applying $J_2$ to $(8)$, we have
\begin{eqnarray}
J^2_2\widetilde{X} & = & J_2GX - \tilde{g} (\widetilde{X}, \xi) \xi - \tilde{g} (\widetilde{X}, J_1\xi) J_1\xi \nonumber \\
&{ }& - g(X, J_2\xi) J_2\xi - g(X, J_3\xi) J_3\xi \nonumber \\
& = & G(GX) + \tilde{g} (J_2 (GX), \xi) \xi + \tilde{g} (J_2 (GX), J_1\xi) J_1\xi - \tilde{g} (\widetilde{X}, \xi) \xi \nonumber \\
&{ }& - \tilde{g} (\widetilde{X}, J_1\xi) J_1\xi - g(X, J_2\xi) J_2\xi - g(X, J_3\xi) J_3\xi \nonumber \\
& = & G^2X - g(GX, J_2 \xi) \xi + g(GX, J_3\xi) J_1\xi - \tilde{g} (\widetilde{X}, \xi) \xi \\
&{ }& - \tilde{g} (\widetilde{X}, J_1\xi) J_1\xi - g(X, J_2\xi) J_2\xi - g(X, J_3\xi) J_3\xi. \nonumber
\end{eqnarray}
On the other hand, by the definition of $J_2$,
\begin{eqnarray}
J^2_2 \widetilde{X} = - \widetilde{X} = - X - \tilde{g} (\widetilde{X}, \xi) \xi - \tilde{g} (\widetilde{X}, J_1\xi) J_1\xi.
\end{eqnarray}
Comparing the tangent and normal components in $(9)$ and $(10)$, we get
\begin{eqnarray}
&{ }& G^2X - g(X, J_2\xi) J_2\xi - g(X, J_3\xi) J_3\xi  = - X, \\
&{ }& g(GX, J_2 \xi) = g(GX, J_3\xi) = 0.
\end{eqnarray}
Now we define 1-forms $u$ and $v$, and unit dual vector fields $U$ and $V$ with respect to $g$ by 
\begin{eqnarray}
&{ }& u(X) = g(X, J_2 \xi), \ \ v(X) = g(X, J_3 \xi) = - (u \circ J)(X), \\ 
&{ }& U = J_2 \xi, \ \ V = J_3\xi = JU.
\end{eqnarray}
By this definitions, $(11)$ and $(12)$ show respectively
\begin{eqnarray}
G^2 = -id. + u \otimes U + v \otimes V, \ \ u(GX) = v(GX) = 0.
\end{eqnarray}
Also, applying $X = U$ and $X = V$ to $(13)$, we have $GU = 0$ and $GV = 0$ respectively. Since $J_2$ is skew-symmetric with respect to $\tilde{g}$, i.e. $\tilde{g}(J_2 X, Y) = - \tilde{g} (X, J_2Y)$, we have $\tilde{g}(GX, Y) = - \tilde{g} (X, GY)$.
Similarly, applying $J_3$ to $(6)$, we have
\begin{eqnarray}
J_3\widetilde{X} & = & HX + \tilde{g} (X, \xi) J_3\xi + \tilde{g} (J_3X, J_1\xi) J_1\xi + \tilde{g}(\widetilde{X}, \xi) J_3\xi + g(\widetilde{X}, J_1\xi) J_2\xi \nonumber \\
& = & HX + g (X, \xi) J_3\xi + g(X, J_2\xi) J_1\xi + \tilde{g}(\widetilde{X}, \xi) J_3\xi + g(\widetilde{X}, J_1\xi) J_2\xi, \nonumber 
\end{eqnarray}
where $HX$ is the component of $J_3X$ tangent to $TM$, and some relations similar to $(12), (13)$ and $(14)$, 
\begin{eqnarray}
&{ }& H^2 = -id. + u \otimes U + v \otimes V, \ \ HU = HV = 0, \\
&{ }& u\circ H = v \circ H = 0, \ \ g(HX, Y) = - g(X, HY).
\end{eqnarray}
Now applying $J_2$ to $(15)$, we have
\begin{eqnarray}
J_2J_3\widetilde{X} & = & J_2HX + \tilde{g} (\widetilde{X}, \xi) J_1\xi - \tilde{g} (\widetilde{X}, J_1\xi) \xi \nonumber \\
&{ }& - g(X, J_3\xi) J_2\xi + g(X, J_2\xi) J_3\xi \nonumber \\
& = & G(HX) - g(HX, J_2\xi) \xi + g(HX, J_3\xi) J_1\xi + \tilde{g} (\widetilde{X}, \xi) J_1\xi \nonumber \\
&{ }& - \tilde{g} (\widetilde{X}, J_1\xi) \xi - v(X)U + u(X)V \nonumber \\
& = & GHX - v(X)U + u(X)V + \tilde{g} (\widetilde{X}, \xi) J_1\xi - \tilde{g} (\widetilde{X}, J_1\xi) \xi.
\end{eqnarray}
On the other hand, 
\begin{eqnarray}
J_2J_3\widetilde{X} = J_1\widetilde{X} = JX + + \tilde{g} (\widetilde{X}, \xi) J_1\xi - \tilde{g} (\widetilde{X}, J_1\xi) \xi. \nonumber
\end{eqnarray}
By comparing the tangent parts in $(17)$ and $(18)$, we get 
\begin{eqnarray}
GH = J + v \otimes U - u \otimes V.
\end{eqnarray}
Finally, applying $G$ from left side to $(14)$ and $(21)$, we have respectively
\begin{eqnarray}
&{ }& G^2H = G(GH) = GJ + v \otimes GU - u \otimes GV = GJ, \nonumber \\
&{ }& G^2H = - H + (u \circ H) \otimes U + (v \circ H) \otimes V = - H, \nonumber
\end{eqnarray}
which show 
\begin{eqnarray}
GJ = - H. 
\end{eqnarray}
From $(13), (14), (15), (16), (17), (19)$ and $(20)$, the structure $(G, H, J, u, v, U, V, g)$ satisfies the definition of complex almost contact metric structure. Then we conclude our theorem. \\ \\
$\bf{Theorem \ 3.3.}$ \  Let $(\widetilde{M}, J_1, J_2, J_3, \tilde{g})$ be a hyperk\"{a}hler manifold and $M$ be a complex hypersurface of $\widetilde{M}$. The inclusion $\iota : M \longrightarrow \widetilde{M}$ canonically induces a complex almost contact metric structure on $M$. \\
\section{connections}

This section is based on a paper by B. Smyth $\cite{S}$. Let $\widetilde{\nabla}$ be the Levi-Civita connection with respect to $\tilde{g}$. For $X, Y \in TM$, we define a tensor field of type $(1,1)$ $A$ and a 1-form $s$ by
\begin{eqnarray*}
\widetilde{\nabla}_X \xi = - AX + s(X) J_1\xi.
\end{eqnarray*}
Then we decompose $\widetilde{\nabla}_X{Y}$ to 
\begin{eqnarray*}
\widetilde{\nabla}_X Y = \nabla_X Y + g(AX, Y) \xi + g(JAX, Y) J_1\xi,
\end{eqnarray*}
where $\nabla_X Y$ denotes the component of $\widetilde{\nabla}_X Y$ tangent to $M$. It is known that $\nabla$ is the Levi-Civita connection with respect to $g$. Now we give expressions for the covariant derivatives of $G$ and $H$ on a complex almost contact metric structure on $M$. \\ \\
$\bf{Theorem \ 4.1}$ \  Let $(G, H, J, u, v, U, V, g)$ be a complex almost contact metric structure on a complex hypersurface $M$. Then the derivatives of $G$ and $H$ have the following forms: \\ \\
\begin{eqnarray*}
&{ }& (\nabla_X G) Y = - u(Y)AX + v(Y)JAX + g(AX, Y)U - g(JAX, Y)V, \\
&{ }& (\nabla_X H) Y = - u(Y)JAX - v(Y)AX + g(AX, Y)V + g(JAX, Y)U. \\
\end{eqnarray*}
$Proof.$ Since from $(7)$, $J_2 Y = GY + \tilde{g} (J_2 Y, \xi) \xi + \tilde{g} (J_2 Y, J_1\xi) J_1\xi$, we have
\begin{eqnarray}
\widetilde{\nabla}_X (J_2Y) & = & \widetilde{\nabla}_X (GY) + \tilde{g} (J_2 \widetilde{\nabla}_X Y, \xi) \xi + \tilde{g} (J_2 Y,  \widetilde{\nabla}_X\xi) \xi + \tilde{g} (J_2 Y,  \xi) \widetilde{\nabla}_X \xi  \nonumber \\
& { } & + \ \tilde{g} (J_2  \widetilde{\nabla}_X Y, J_1\xi) J_1\xi + \tilde{g} (J_2 Y, J_1 \widetilde{\nabla}_X \xi) J_1\xi + \tilde{g} (J_2 Y, J_1\xi) J_1 \widetilde{\nabla}_X \xi \nonumber \\
& = & \nabla_X (GY) + g(AX, GY) \xi + g(JAX, GY) J_1\xi - \tilde{g}(\widetilde{\nabla}_X Y, J_2 \xi) \xi \nonumber \\
& { } & + \ \tilde{g} (J_2 Y, -AX + s(X)J_1\xi) \xi - \tilde{g} (Y, J_2 \xi)(- AX + s(X)J_1 \xi) \nonumber \\
& { } & + \ \tilde{g} (\widetilde{\nabla}_X Y, J_3 \xi) J_1\xi + \tilde{g} (J_2 Y, - JAX - s(X)\xi) J_1\xi \nonumber \\
& { } & + \ \tilde{g} (Y, J_3\xi) (- JAX - s(X)J_1\xi) \nonumber \\
& = &  \nabla_X (GY) + g(AX, GY) \xi + g(JAX, GY) J_1\xi - u(\nabla_X Y) \xi \nonumber \\
& { } & - \ g(GY, AX) \xi + s(X)v(Y) \xi + u(Y)AX - s(X)u(Y)J_1\xi \nonumber \\
& { } & + \ v(\nabla_X Y) J_1\xi - g(GY, JAX) J_1\xi + s(X)u(Y) J_1\xi  \nonumber \\
& { } & - \ v(Y) JAX - s(X)v(Y) \xi \nonumber \\
& = & (\nabla_X G)Y + G \nabla_X Y - u(\nabla_X Y) \xi + u(Y) AX + v(\nabla_X Y) J_1 \xi \\
& { } & - \ v(Y) JAX. \nonumber 
\end{eqnarray}
On the other hand, since $J_2$ is parallel to $\widetilde{\nabla}$, we have
\begin{eqnarray}
\widetilde{\nabla}_X (J_2Y) & = & J_2 \widetilde{\nabla}_X Y \nonumber \\
& = & J_2 (\nabla_X Y + g(AX, Y) \xi + g(JAX, Y) J_1\xi) \nonumber \\
& = & G\nabla_X Y + \tilde{g} (J_2 \nabla_X Y, \xi) \xi + \tilde{g} (J_2 \nabla_X Y, J_1\xi) J_1\xi \nonumber \\
& { } & + \ g(AX, Y) \xi - g(JAX, Y) J_1\xi \nonumber \\
& = & G\nabla_X Y - u(\nabla_X Y) \xi + v(\nabla_X Y) J_1\xi +  g(AX, Y) \xi - g(JAX, Y) J_1\xi
\end{eqnarray}
Comparing $(21)$ and $(22)$, we get
\begin{eqnarray}
(\nabla_X G) Y & = & - u(Y)AX + v(Y)JAX \\
&{ }& + g(AX, Y)U - g(JAX, Y)V. \nonumber
\end{eqnarray}
Since $J$ is parallel to $\nabla$, we have
\begin{eqnarray*}
J (\nabla_X G) Y = \nabla_X (JGY) - JG \nabla_X Y = (\nabla_X H)Y.
\end{eqnarray*}
Then applying $J$ to $(23)$, we get
\begin{eqnarray*}
(\nabla_X H) Y = - u(Y)JAX - v(Y)AX + g(AX, Y)V + g(JAX, Y)U. \\
\end{eqnarray*}
$\bf{Proposition \ 4.2}$ \  Let $(G, H, J, u, v, U, V, g)$ be a complex almost contact metric structure on a complex hypersurface $M$. Then 
\begin{eqnarray*}
s(X) = g(\nabla_X V, U) = - \sigma(X). \\
\end{eqnarray*}
$Proof.$ \ By Proposition 4.1, we have
\begin{eqnarray}
G(\nabla_X G) Y & = & \nabla_X (G^2 Y) - G^2\nabla_X Y - (\nabla_X G)(GY) \nonumber \\
& = & \nabla_X (- Y + u(Y)U + v(Y)V) + \nabla_X Y - u(\nabla_X Y)U - v(\nabla_X Y)V \nonumber \\
& { } & - \ g(AX, GY)U + g(JAX, GY)V \nonumber \\
& = & - \nabla_X Y + g(\nabla_X Y, U)U + g(Y, \nabla_X U)U + g(Y, U)\nabla_X U \nonumber \\
& { } & + \ g(\nabla_X Y, V)V + g(Y, \nabla_X V)V + g(Y, V) \nabla_X V + \nabla_X Y \nonumber \\
& { } & - \ u(\nabla_X Y)U - v(\nabla_X Y)V - g(AX, GY)U + g(JAX, GY)V \nonumber \\
& = & g(Y, \nabla_X U)U + u(Y)\nabla_X U + g(\nabla_X V, Y)V + v(Y) \nabla_X V \\
& { } & + \ g(GAX, Y) U + g(HAX, Y) V. \nonumber 
\end{eqnarray}
On the other hand, by proposition 4.1, we have 
\begin{eqnarray}
G(\nabla_X G) Y = - u(Y) GAX - v(Y)HAX. 
\end{eqnarray}
Comparing $(24)$ and $(25)$, we get 
\begin{eqnarray}
& { } & g(Y, \nabla_X U)U + u(Y)\nabla_X U + g(\nabla_X V, Y)V + v(Y) \nabla_X V \\
& { } & + \ g(GAX, Y) U + g(HAX, Y) V + u(Y) GAX - v(Y)HAX = 0. \nonumber
\end{eqnarray}
Applying $u$ to $(26)$, we have $\nabla_X U = - GAX - u(\nabla_X V)V$, since for $X$ and $Y$,
\begin{eqnarray*}
0 & = & g(Y, \nabla_X U) + v(Y) u(\nabla_X V) + g(GAX, Y) \\
& = & g(Y, \nabla_X U + u(\nabla_X V)V + GAX).
\end{eqnarray*}
By this expression, we get
\begin{eqnarray}
\widetilde{\nabla}_X U & = & \nabla_X U + g(AX, U) \xi + g(JAX, U) J_1 \xi \\
& = & - \ GAX - u(\nabla_X V)V + u(AX) \xi - v(AX) J_1\xi. \nonumber\end{eqnarray}
On the other hand, since $J_2$ is parallel to $\widetilde{\nabla}$, we have
\begin{eqnarray}
\widetilde{\nabla}_X U & = & \widetilde{\nabla}_X (J_2 \xi) \nonumber \\
& = & J_2 ( - AX + s(X) J_1 \xi ) \nonumber \\
& = & - \ GAX - \tilde{g} (J_2AX, \xi) \xi + \tilde{g} (J_2AX, J_1\xi) J_1 \xi - s(X) J_3 \xi \nonumber \\
& = & - \ GAX + u(AX) \xi - v(AX) J_1\xi - s(X) V.
\end{eqnarray}
Comparing $(27)$ and $(28)$, we get the conclusion. \\ \\
$\bf{Proposition \ 4.3}$ \  For any $X \in TM$, $\nabla_X G$ and $\nabla_X H$ are skew-symmetric operators with respect to $g$. \\ \\
$Proof.$ \ By proposition 4.1, we get the following equality, which gives the conclusion.
\begin{eqnarray*}
g( (\nabla_X G) Y, Z) & = & - \ u(Y)g(AX, Z) + v(Y)g(JAX, Z) \\
& { } & + \ g(AX, Y)u(Z) - g(JAX, Y)v(Z) \\
& = & - \ g((\nabla_X G) Z, Y). 
\end{eqnarray*}
Department of Natural Sciences \\ 
National Institute of Technology, Ibaraki Collage \\ 
886 Nakane, Hitachinaka, Ibaraki 312-8508 Japan \\ 
Tel: +81-29-271-2865
e-mail: imadam@ge.ibaraki-ct.ac.jp
\end{document}